\documentclass[a4paper, 10pt]{article}

\usepackage{amsmath} 
\usepackage{amsfonts} 
\usepackage{amsthm} 
\usepackage{epsf}  
\usepackage{appendix}
\usepackage{listings}
\usepackage{url}

\newtheorem{thm}{Theorem}

\theoremstyle{definition}
\newtheorem{defn}[thm]{Definition}

\theoremstyle{remark}
\newtheorem*{rem}{Remark}

\def\address#1#2{\begingroup
\noindent\parbox[t]{7.8cm}{%
\small{\scshape\ignorespaces#1}\par\vskip1ex
\noindent\small{\itshape E-mail address}%
\/: #2\par\vskip4ex}\hfill%
\endgroup}%

\title{Hilbert polynomial of the Kimura 3-parameter model}
\author{Kaie Kubjas}
\date{}

\begin{document}

\newcommand{\T}[1]{{#1}_{T}}
\newcommand{\cat}[1]{{#1}_{\text{3c}}}
\newcommand{\sn}[1]{{#1}_{\text{sn}}}
\newcommand{\three}[1]{{#1}_{\text{3l}}}
\newcommand{\four}[1]{{#1}_{\text{4l}}}
\newcommand{\ncat}[2]{{#1}_{\text{3c}}^{#2}}
\newcommand{\nsn}[2]{{#1}_{\text{sn}}^{#2}}
\newcommand{\nthree}[2]{{#1}_{\text{3l}}^{#2}}
\newcommand{\nfour}[2]{{#1}_{\text{4l}}^{#2}}

\maketitle

\begin{abstract}
In \cite{BW07} Buczy\'{n}ska and Wi\'{s}niewski showed that for the Jukes Cantor binary model of a 3-valent tree the Hilbert polynomial depends only on the number of leaves of the tree and not on its shape. We ask if this can be generalized to other group-based models. The Jukes Cantor binary model has $\mathbb{Z}_2$ as the underlying group. We consider the Kimura 3-parameter model with $\mathbb{Z}_2 \times \mathbb{Z}_2$ as the underlying group. We show that the generalization of the statement about the Hilbert polynomials to the Kimura 3-parameter model is not possible as the Hilbert polynomial depends on the shape of a 3-valent tree.
\end{abstract}

\section{Introduction}
Phylogenetic algebraic geometry studies complex algebraic varieties arising from evolutionary models in biology (see for example \cite{ESSR05}, \cite{SS05}). The Jukes Cantor binary model is the simplest model. Buczy\'{n}ska and Wi\'{s}niewski showed in \cite{BW07} that the Hilbert polynomial of the ideal of the Jukes Cantor binary model with a 3-valent tree $T$ depends only on the number of leaves of $T$ and not on the shape. We are interested if this property of the Hilbert polynomial can be generalized to more complicated models.

The Jukes Cantor binary model is a group-based model with the underlying group $\mathbb{Z}_2$. The most natural generalization seems to be the Kimura 3-parameter model, which is a group-based model with the underlying group \mbox{$\mathbb{Z}_2 \times \mathbb{Z}_2$}. However, in this paper we conclude that this generalization is not possible. For the Kimura 3-parameter model the Hilbert polynomial also depends on the shape of a 3-valent tree. We show that the ideals of 2 different trees with 8 leaves -- the 3-caterpillar and snowflake trees (see figures on page 3) -- have different Hilbert polynomials. The Kimura 3-parameter model being the closest model to the DNA binary model, it is unlikely that the property about Hilbert polynomials would hold for other models.

In Section 2 we recall the construction for the Kimura 3-parameter model. In Section 3 we show that the Hilbert polynomials of the ideals of the \mbox{3-caterpillar} and snowflake trees have different values when evaluated at 3 and hence the Hilbert polynomials are not the same. The main idea is to decompose the original trees to smaller trees and use toric fiber products introduced by Sullivant in \cite{Sullivant07}. Michalek showed in \cite{Michalek10} that the Kimura 3-parameter model is normal, so we reduce the problem of evaluating the Hilbert polynomials of toric varieties to evaluating the Ehrhart polynomials of the corresponding polytopes. Computations are done with the help of \texttt{polymake} \cite{polymake}, \cite{JMP09} and \texttt{Normaliz} \cite{normaliz}.

I would like to thank Christian Haase and Andreas Paffenholz for their ideas and technical help.

\section{Kimura 3-parameter model}
We will not give the parametric construction of the Kimura 3-parameter model coming from biology here (see \cite{SS05}), but will directly define the toric ideal and the corresponding polytope of the Kimura 3-parameter model. By doing this we follow \cite{SS05} and \cite{Sullivant07}. 

Let $T$ be a 3-valent tree with $n+1$ leaves labeled by $1,\ldots, n+1$, let the root be at the leaf $n+1$, and direct the edges away from the root. A leaf $l$ is a descendant of an edge $e$ if there is a directed path from $e$ to $l$. Denote by de($e$) the set of all descendants of the edge $e$.

For a sequence $g_1,\ldots ,g_n$ in $\mathbb{Z}_2 \times \mathbb{Z}_2$, we define
\begin{displaymath}
g_e=\sum_{i\in \text{de}(e)}g_i,
\end{displaymath}
where $e$ is an edge of $T$ and the subindices $i$ denote simultaneously leaves and their labels. Let
\begin{displaymath}
\mathbb{K}[q]=\mathbb{K}[q_{g_1,\ldots ,g_n}|g_i\in \mathbb{Z}_2 \times \mathbb{Z}_2] \text{ and } \mathbb{K}[a]=\mathbb{K}[a_h^{(e)}|e\in E(T),h\in \mathbb{Z}_2 \times \mathbb{Z}_2]
\end{displaymath}
and consider the ring homomorphism
\begin{align*}
\phi _T:\mathbb{K}[q] & \rightarrow \mathbb{K}[a]\\
q_{g_1,\ldots ,g_n} & \mapsto \prod _{e\in E(T)} a_{g_e}^{(e)}.
\end{align*}

\begin{defn}
Let the ideal $\T{I}=\text{ker}(\phi _T)$ be the ideal of the Kimura 3-parameter model with tree $T$.
\end{defn}

The ideal $\T{I}$ is a toric ideal and we can define the corresponding polytope.

\begin{defn}
Let the polytope 
\begin{displaymath}
\T{P}=\text{conv} (\{\alpha \in \mathbb{Z}^{E(\mathcal{T})\times (\mathbb{Z}_2 \times \mathbb{Z}_2)}|a^{\alpha }=\phi _T(q_{g_1,\ldots ,g_n}), q_{g_1,\ldots ,g_n}\in \mathbb{K}[q]\})
\end{displaymath} 
be the polytope of the Kimura 3-parameter model with tree $T$.
\end{defn}

\begin{defn}
Let the lattice $\T{L} \subseteq \mathbb{Z}^{E(\mathcal{T})\times (\mathbb{Z}_2 \times \mathbb{Z}_2)}$ be the lattice generated by the vertices of $\T{P}$.
\end{defn}

Since $T$ is an acyclic directed graph, there is an induced partial order on the edges of $T$. Namely $e<e'$ if there is a directed path from $e'$ to $e$. Let $T$ be a tree that contains an interior edge $e$. Then $e$ induces a decomposition of $T$ as $T_e^+*T_e^-$ where $T_e^-$ is a subtree of $T$ consisting of all edges $e'\in T$ with $e'\leq e$ and $T_e^+$ consists of all edges $e'\in T$ with $e'\not < e$. Thus $T_e^+$ and $T_e^-$ overlap in the single edge $e$. We root $T_e^-$ by the tail of $e$, and keep the root of $T_e^+$ at the original root $n+1$. Let $I_{T^+}^e$ and  $I_{T^-}^e$ denote the ideals of the Kimura 3-parameter model with trees $T_e^+$ and $T_e^-$ and let $\mathbb{K}[q]_+$ and $\mathbb{K}[q]_-$ denote the ambient polynomial rings, respectively.

Using toric fiber products of ideals (see \cite{Sullivant07}) Sullivant stated the following theorem.

\begin{thm}[Sullivant \cite{Sullivant07}]
Let $T$ be a tree with an interior edge $e$, and resulting decomposition $T=T_e^+*T_e^-$. For each variable $q_{\mathbf{g}}=q_{g_1,\ldots ,g_n}$ in $\mathcal{K}[q], \mathcal{K}[q]_+$ and $\mathcal{K}[q]_-$, let deg$(q_{\mathbf{g}})=e_{g_e}$. Then 
\begin{displaymath}
\T{I}=I_{T^+}^e\times _{\mathcal{A}}I_{T^-}^e
\end{displaymath}
with $\mathcal{A}=\{e_{(0,0)},e_{(0,1)},e_{(1,0)},e_{(1,1)}\}$.
\end{thm}

\section{Counting lattice points}

\begin{enumerate}
\item Since the polytopes of the 3-caterpillar and snowflake trees are too large to compute their lattice points directly, 
we decompose them into smaller trees like shown on the figure below.
\begin{center}
\epsfbox{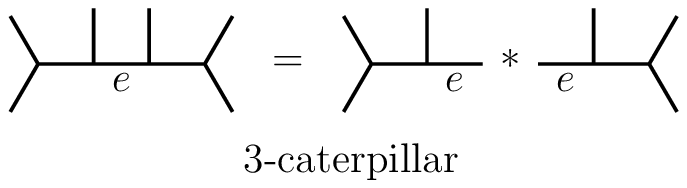}

\epsfbox{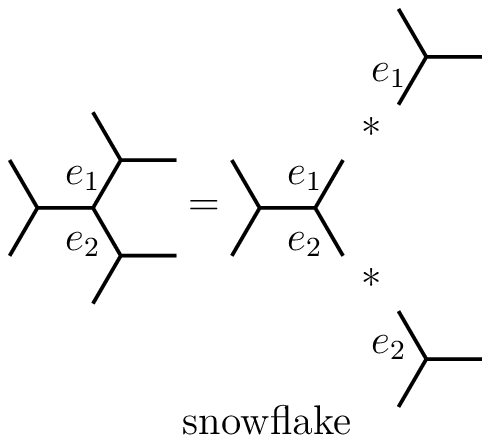}
\end{center}

Henceforth we use the abbreviations 3c, sn, 3l, 4l for the 3-caterpillar, snowflake, 3-leaf and 4-leaf trees, respectively.

In the decomposition of the 3-caterpillar tree define deg$(q_{\mathbf{g}})=e_{g_e}$ for $q_{\mathbf{g}}$ in $\cat{\mathcal{K}[q]}$ and  $\four{\mathcal{K}[q]}$. Then
\begin{gather*}
\cat{I}=\nfour{I}{e}\times _{\mathcal{A}} \nfour{I}{e}
\end{gather*}
with $\mathcal{A}=\{e_{(0,0)},e_{(0,1)},e_{(1,0)},e_{(1,1)}\}$.

In the decomposition of the snowflake tree define deg($q_{\mathbf{g}})=e_{g_{e_1},g_{e_2}}$ for $q_{\mathbf{g}}$ in  $\sn{\mathcal{K}[q]}$ and $\four{\mathcal{K}[q]}$ and define deg($q_{\mathbf{g}})=e_{g_{e_i}}$ with $i\in \{ 1,2 \}$ for $q_{\mathbf{g}}$ in $\three{\mathcal{K}[q]}$. Then
\begin{gather*}
\sn{I}=\nfour{I}{e_1,e_2}\times _{\mathcal{A}}\nthree{I}{e_1}\times _{\mathcal{A}}\nthree{I}{e_2}
\end{gather*}
with $\mathcal{A}=\{e_{(0,0)},e_{(0,1)},e_{(1,0)},e_{(1,1)}\}$.

\item Denote the multigraded Hilbert function of $\mathbb{K}[q]/I$ by $h(\mathbb{K}[q]/I;u)$. In \cite{Sullivant07} Sullivant gives a formula for computing multigraded Hilbert functions of toric fiber products. Applying this to the decompositions of Step 1 gives for $u,v\in \mathbb{Z}^{\mathbb{Z}_2\times \mathbb{Z}_2}$
\begin{gather*}
h(\cat{\mathbb{K}[q]}/\cat{I};u)
=h(\four{\mathbb{K}[q]}/\nfour{I}{e};u)h(\four{\mathbb{K}[q]}/\nfour{I}{e};u),\\
h(\sn{\mathbb{K}[q]}/\sn{I};u,v)
=h(\four{\mathbb{K}[q]}/\nfour{I}{e_1,e_2};u,v)h(\three{\mathbb{K}[q]}/\nthree{I}{e_1};u)h(\three{\mathbb{K}[q]}/\nthree{I}{e_2};v).
\end{gather*}

\item A monomial having multidegree $u\in \mathbb{Z}^{\mathbb{Z}_2\times \mathbb{Z}_2}$ has total degree $\sum _{h\in \mathbb{Z}_2 \times \mathbb{Z}_2}u_h$. Thus single graded Hilbert functions can be computed using multigraded Hilbert functions
\begin{gather*}
h(\cat{\mathbb{K}[q]}/\cat{I};n)=\sum _{u:\sum u_h=n}h(\cat{\mathbb{K}[q]}/\cat{I};u),\\
h(\sn{\mathbb{K}[q]}/\sn{I};n)=\sum _{u,v:\sum u_h=n,\sum v_h=n}h(\sn{\mathbb{K}[q]}/\sn{I};u,v).
\end{gather*}

\item In \cite{Michalek10} Michalek shows that the polytopes of the Kimura 3-parameter model polytopes are normal, hence corresponding Ehrhart and Hilbert polynomials are equal (see for example \cite{Sturmfels96}). Thus $h(\mathbb{K}[q]/I_T;u)$ counts lattice points in the lattice $\T{L}$ of the $\sum _{h\in \mathbb{Z}_2 \times \mathbb{Z}_2}u_h$ dilation of the polytope $\T{P}$ intersected with hyperplanes $\{x_h^e=u_h\}, h\in \mathbb{Z}_2 \times \mathbb{Z}_2$. Using Step 2 and Step 3 we get
\begin{eqnarray*}
\text{ehr}_{\cat{P}}(n)&=&\sum _{u:\sum u_h=n} \left\vert n\four{P}\bigcap\{x_h^e=u_h\}\bigcap \four{L}||n\four{P}\bigcap \{x_h^e=u_h\}\bigcap \four{L}\right\vert, \\
\text{ehr}_{\sn{P}}(n)&=&\sum _{u,v:\sum u_h=n,\sum v_h=n} \left\vert n\four{P}\bigcap \{ x_h^{e_1}=u_h \} \bigcap \{x_h^{e_2}=v_h\}\bigcap \four{L} \right\vert \\
& \cdot & \left\vert n\three{P}\bigcap \{ x_h^{e_1}=u_h\}\bigcap \three{L}\right\vert \left\vert n\three{P}\bigcap \{x_h^{e_2}=v_h\}\bigcap \three{L}\right\vert.
\end{eqnarray*}

\item Using \texttt{polymake} and \texttt{Normaliz} we can compute $|3\T{P}\cap \{x_l^e=u_l\} \cap \T{L}|$ for 3-leaf and 4-leaf trees. It is important to do the basis transformation before counting lattice points, since these programs assume that the lattice is the standard lattice. Using formulas from Step 4 we got that in the 3rd dilation the 3-caterpillar polytope has 69324800 and the snowflake polytope has 69248000 lattice points. Hence their Ehrhart (and thus Hilbert) polynomials are different.

\begin{rem}
Similar computations show that in the 2nd dilation the 3-caterpillar and snowflake polytopes have both 396928 lattice points.
\end{rem}

\end{enumerate}

\bibliographystyle{alpha}
\bibliography{K3par}

\address{Mathematisches Institut\\
Freie Universit\"at Berlin\\
Arnimallee 3\\
14195 Berlin, Germany}{kubjas@math.fu-berlin.de}

\end{document}